\documentclass[12pt]{amsart}

\usepackage{amsfonts,amsmath,amsthm}
\usepackage{latexsym}

\newtheorem{theorem}{Theorem}[section]
\newtheorem{corollary}[theorem]{Corollary}

\theoremstyle{definition}

\newcommand \re {\mathrm{Re}\:}
\newcommand\kr{\mathrm{Ker}\:}
\newcommand\dom{\mathrm{Dom}\:}

\newcommand{\R}{\ensuremath{\mathbb{R}}}

\def \d {\delta}
\def \l {\lambda}

\def \f {\phi}

\def \< {\langle}
\def \> {\rangle}

\begin{document}

\title[On quasi-contractivity of $C_0$-semigroups on Banach spaces]{On quasi-contractivity of $C_0$-semigroups on Banach spaces}
\author{M\'at\'e~Matolcsi}
\email{matomate@renyi.hu}
\address{ Alfr\'ed R\'enyi Institute of Mathematics,
                                     Hungarian Academy of Sciences
                                     POB 127
                                     H-1364 Budapest, Hungary
                                     Tel: (+361) 483-8302
                                     Fax: (+361) 483-8333}

\date{\today}
\maketitle

\begin{abstract}
A basic result in semigroup theory states that every $C_0$-semigroup is quasi-contractive with respect to some appropriately chosen equivalent norm. This paper contains a counterpart of this well-known fact. Namely, by examining the convergence of the Trotter-type formula $(e^{\frac{t}{n}A}P)^n$ (where $P$ denotes a bounded projection), we prove that whenever the generator $A$ is unbounded it is possible to introduce an equivalent norm on the space with respect to which the semigroup is {\it{not}} quasi-contractive.

\end{abstract}

\section{Introduction}
Many important results in semigroup theory  rely on the simple fact that for a given $C_0$-semigroup $e^{tA}$ on a Banach space $X$ it is always possible to introduce an equivalent norm on $X$ with respect to which $e^{tA}$ is quasi-contractive. While examining the convergence of the  Trotter-type formula 
\begin{equation}\label{pro}
(e^{\frac{t}{n}A}P)^n
\end{equation}
(see \cite{mat}), the author was led to the natural question of whether it is
always possible to introduce an equivalent norm on $X$ with respect to which
$e^{tA}$ is {\it{not}} quasi-contractive. This is clearly not possible if the
generator $A$ is bounded. However, if  $A$ is unbounded then it
is natural to expect that such a norm does exist. Indeed, in \cite{mat},
Theorem 2, the Hilbert space version of this question was settled: assuming
that $A$ is unbounded an equivalent {\it{scalar product}} (not merely a norm)
was constructed with respect to which $e^{tA}$ is not quasi contractive. This
result was then used to show that whenever $A$ is unbounded it is
possible to find a bounded (but not necessarily orthogonal) projection $P$
such that \eqref{pro} does not converge strongly (cf. \cite{mat}, first part
of Theorem 3). The proofs of these results, however, relied heavily on the
notion of orthogonality. The aim of this paper is to prove the Banach space
analouge of these two results (see Theorem \ref{thm1} and Corollary \ref{cor1}
below). 

We are aware that the existence of a 'non-quasi-contractive' norm will
probably not have as much use as that of a 'quasi-contractive'
one. Nevertheless, it gives an affirmative answer to a natural
question, and shows that whenever $A$ is unbounded it is up to our choice
whether to regard $e^{tA}$ as being quasi-contractive or non-quasi-contractive.

As the motivation to tackle the questions above came from the investigations
of the convergence of formula \eqref{pro}, we mention that
the history of this formula goes back to \cite{valaki}, \cite{dav},
\cite{kato} and \cite{ab}. The interested reader can also find a brief
overview and some recent results in \cite{ms} and \cite{mat}. Thus far, most
of the attention concerning formula \eqref{pro} has been devoted to the
Hilbert space case, but Theorem \ref{thm1} below shows that some of the results remain true in the most general settings. 

\section{Main results}

In the Banach space setting we reverse the steps of \cite{mat}. First we
characterize the convergence of \eqref{pro} in terms of the generator $A$, and
then we use this result to construct an equivalent norm with respect to which
the semigroup $e^{tA}$ is not quasi-contractive. 

\begin{theorem}\label{thm1}
Let $A$ generate a $C_0$-semigroup $e^{tA}$ on a complex Banach space $X$. The following are equivalent: 

(i) $A$ is bounded 

(ii) $\lim_{n\to \infty}(e^{\frac{t}{n}A}P)^nx$ exists for all bounded projection $P$, and all $x\in X$, $t\ge 0$.

\begin{proof}

The implication $(i)\rightarrow (ii)$ is fairly standard and contained in
\cite{ms}, Theorem 1.  

For the implication $(ii)\rightarrow (i)$ assume that $A$ is unbounded. Then,
there exists an element $\phi\in X^\ast$ such that  $\phi\notin \dom
A^\ast$. $\kr \phi$, the kernel of $\phi$, is a 1-codimensional subspace of
$X$. As $\dom A$ is dense in $X$, $\dom A \not\subset \kr\phi$. Therefore,  we
can choose a vector $x\in \dom A$ such that $\phi (x)=1$. Let $P_x$ denote the  projection along $\kr \phi$ onto the 1-dimensional subspace spanned by $x$; i.e. we decompose each element $z\in X$ as $z=\f (z)x +(z-\f (z)x)$ and we let $ P_xz=\f (z)x$. Note that $\f (P_x z)=\f (z)$. Therefore, 
\begin{equation}\label{eq2}
 \f ((e^{\frac{1}{n}A}P_x)^nx)= \f ((P_xe^{\frac{1}{n}A}P_x)^nx)
\end{equation}

Now, observe that $(P_xe^{\frac{1}{n}A}P_x)x=c_nx$ where $c_n=\f (e^{\frac{1}{n}A}x)$. Therefore, 
\begin{equation}\label{eq3}
\f ((P_xe^{\frac{1}{n}A}P_x)^nx)= \f (c_n^nx)=c_n^n=(\f (e^{\frac{1}{n}A}x))^n
\end{equation}

Furthermore, 
$$ n(c_n-1)= \f\left ( \frac{e^{\frac{1}{n}A}x-x}{1/n}\right )= 
\f\left ( \frac{(e^{\frac{1}{n}A}-I)x}{1/n}\right )$$
therefore 
\begin{equation}\label{eq4}
\lim_{n\to \infty} n(c_n-1)=\f (Ax)
\end{equation}

Combining \eqref{eq2} , \eqref{eq3} and \eqref{eq4} we get
\begin{equation}\label{eq5}
\lim_{n\to\infty}\f\ ((e^{\frac{1}{n}A}P_x)^nx)= \lim_{n\to \infty}c_n^n=e^{\f (Ax)}
\end{equation}
The rest of the proof is similar to the Hilbert space version (see \cite{mat}). 

We are going to construct an element $y\in X$ such that $\f (y)=1$ and $\lim_{n\to\infty}\f\ ((e^{\frac{1}{n}A}P_y)^ny)$ does not exist. The vector $y$ will be given as $\lim_{n\to \infty}x_k$ where the sequence $(x_k)$ is to be constructed in the sequel. 

We would like to choose $x_0$ so that $x_0\in \dom A$, $\f (x_0)=1$ and $\re \f (Ax_0)\ge 0$ holds. To do this we take a vector $z$ such that $z\in \dom (A)$, $\f (z)=1$, and we are looking for $x_0$ in a small neighbourhood of $z$. As $\f \notin \dom A^\ast$ we can find a vector $v\in \dom A$ such that $\|v\| \le \frac{1}{2\|\f\|}$ and $|\f (Av)|\ge 4|\f (Az)|$. Now, let $v':=e^{i\alpha}v$ with suitable $\alpha$ such that $\f (Av')$ is nonnegative real. Let $x_0:=\frac{z+v'}{\f (z+v')}$. Then $x_0\in \dom A$, $\f (x_0)=1$ and $\re \f (Ax_0)=\re (\frac{1}{1+\f (v')}\f (Az))+ \re (\frac{1}{1+\f (v')}\f (Av'))\ge -2|\f (Az)|+ \frac{1}{2}\f (Av')\ge 0$ as desired. 

We know that  $\lim_{n\to\infty}\f\ ((e^{\frac{1}{n}A}P_{x_0})^nx_0)=e^{\f (Ax_0)}$. Let $\epsilon >0$ be fixed. Take an index $n_0$ so large that $|\f ((e^{\frac{1}{n_0}A}P_{x_0})^{n_0}x_0)-e^{\f (Ax_0)}|<\epsilon $.
It is clear from standard continuity arguments that there exists a radius $\delta_0 >0$ such that the conditions $h\in B(x_0,\delta_0)$ and $\f (h)=1$ together imply that $|\f ((e^{\frac{1}{n_0}A}P_{h})^{n_0}h)-e^{\f (Ax_0)}|<2\epsilon $. Without loss of generality we can assume that $\delta_0< \frac{\|x_0\|}{2}$. 

We are going to construct the sequence $(x_k)$ inductively. Assume,  therefore, that vectors $x_0, \ x_1, \dots , x_k$, positive numbers $\d_0, \ \d_1, \dots , \d_k$ and indices $n_0, \ n_1, \dots , n_k$ are already given such that for all $0\le j\le k$ the following hold: 

$x_j\in \dom A$, $\f (x_j)=1$, $\re \f (Ax_j)\ge j$ and  $|\f ((e^{\frac{1}{n_j}A}P_{h})^{n_j}h)-e^{\f (Ax_j)}|<2\epsilon $ for all vectors $h$ satisfying  $h\in B(x_j,\delta_j)$ and $\f (h)=1$.

Assume, furthermore, that $\|x_{j+1}-x_j\|<\min \{\frac{\d_0}{2^{j+1}},  \ \frac{\d_1}{2^j}, \dots \frac{\d_j}{2}\}$ for all $0\le j\le k-1$. 

Clearly, there exists a (sufficiently small) radius $\gamma_k >0$ such that for all $g\in B(0,\gamma_k)$ we have $\frac{x_k+g}{\f (x_k+g)}\in B(x_k, \d )$, where $\d:= \min \{\frac{\d_0}{2^{k+1}},  \ \frac{\d_1}{2^k}, \dots \frac{\d_k}{2}\}$. 
Now, the construction of $x_{k+1}$ from the given vector $x_k$ goes the same
way as
the construction above of $x_0$ from the given vector $z$; using that $\f \notin \dom
A^\ast$ we find an appropriate vector $g\in \dom A$, $\|g\|<\gamma_k$ such
that $\f (Ag)$ is positive and 'large', assuring that the definition
$x_{k+1}:= \frac{x_k+g}{\f (x_k+g)}$ gives $\re \f (Ax_{k+1})\ge k+1$. Note,
also, that $x_{k+1}\in B(x_k,\d)$ because $\|g\|<\gamma_k$. Finally,
the index $n_{k+1}$ and the radius $\d_{k+1}$ are chosen to correspond to the
vector $x_{k+1}$, so that  $|\f
((e^{\frac{1}{n_{k+1}}A}P_{h})^{n_{k+1}}h)-e^{\f (Ax_{k+1})}|<2\epsilon$ holds for all vectors $h$ satisfying  $h\in B(x_{k+1},\delta_{k+1})$ and $\f (h)=1$. 

It is clear, by construction, that the sequence $(x_k)$ converges in $X$. For
$y:=\lim x_k$ we have $\f (y)=1$ (and $\|y\|\ge \frac{\|x_0\|}{2}$ by the
choice of $\d_0$). It is also clear, by construction, that $y\in B(x_k, \d_k)$
for all $k\ge 0$. Hence, for all $k\ge 0$ we have  
$$|\f ((e^{\frac{1}{n_{k}}A}P_{y})^{n_{k}}y)-e^{\f (Ax_{k})}|<2\epsilon $$
Notice that $|e^{\f (Ax_{k})}|=e^{\re \f (Ax_{k})}\ge e^k$. This means that
the sequence $(e^{\frac{1}{n}A}P_{y})^{n}y$ does not converge (even weakly).
\end{proof} 
\end{theorem}

From this result the existence of a 'non-quasi-contractive' norm follows
easily. 

\begin{corollary}\label{cor1}
Let $A$ generate a $C_0$-semigroup $e^{tA}$ on a complex Banach space $X$. The
following are equivalent:

(i) A is bounded

(ii) the semigroup $e^{tA}$ is quasi-contractive with respect to every
equivalent norm on $X$. 

\begin{proof}

The implication $(i)\rightarrow (ii)$ is obvious. 

For the implication $(ii)\rightarrow (i)$ assume that $A$ is not bounded.
By
the proof of the preceding theorem we can find vectors $\f\in X^\ast$, $y\in X$
such that  $\f ((e^{\frac{1}{n}A}P_{y})^{n}y)$ does not converge. Moreover,
for the subsequence $n_k$ above we have $|\f ((e^{\frac{1}{n_k}A}P_{y})^{n_k}y)|\ge
e^k-2\epsilon$. 

Introduce a new norm $\| \ \cdot \ \|_0$ on $X$ by
$\|z\|_0:=\|P_yz\|+\|(I-P_y)z\|$. It is clear that the norms $\| \ \cdot \
\|_0$ and $\| \ \cdot \ \|$ are equivalent, and $P_y$ is contractive with
respect to $\| \ \cdot \ \|_0$. We claim, that $e^{tA}$ is not
quasi-contractive with respect to $\| \ \cdot \ \|_0$. Indeed, assume, by
contradiction, that there exists a $\l \in \R$ such that $\|e^{tA}\|_0\le
e^{\l t}$ for all $t\ge 0$. Then $|\f ((e^{\frac{1}{n_k}A}P_{y})^{n_k}y)|\le
\|\f\|\cdot \|y\|\cdot e^\l$, a contradiction.   
\end{proof}
\end{corollary}

\end{document}